\documentclass{raex}

\usepackage{amssymb,amsthm,amsmath}

\begin{Author} 
	\FirstName{Miroslav}\LastName{Ka\v cena}
	\PostalAddress{Department of Mathematical Analysis, Faculty of Mathematics and Physics, Charles University, Sokolovsk\'a~83, 186~75 Praha~8, Czech Republic}
	\Email{kacena@karlin.mff.cuni.cz }
	\Thanks{The author was supported by the grant GAUK 126509.}
\end{Author}

\begin{Author} 
	\FirstName{Luca}\LastName{Motto Ros}
	\PostalAddress{Albert-Ludwigs-Universit\"at Freiburg,
Mathematisches Institut, Abteilung f\"ur Mathematische Logik,
Eckerstra{\ss}e 1,
 D-79104 Freiburg im Breisgau,
Germany}
	\Email{luca.motto.ros@math.uni-freiburg.de}
\end{Author}

\begin{Author} 
	\FirstName{Brian}\LastName{Semmes}
	\PostalAddress{(formerly) Institute for Logic, Language and Computation,
University of Amsterdam, Holland}
	\Email{bsemmes@ovi.com}
\end{Author}

\begin{MathReviews} 
	\primary{03E15} 
	\primary{54H05}
\end{MathReviews}

\begin{KeyWords} 
\end{KeyWords}

\title{Some observations on `A new proof of a theorem of Jayne and Rogers'}

\newcommand{\R}{\RealNumbers}
\newcommand{\Q}{\RationalNumbers}
\newcommand{\N}{\NaturalNumbers}

\theoremstyle{plain}
\newtheorem{thm}{Theorem}
\newtheorem*{thm*}{Theorem}
\newtheorem{lem}[thm]{Lemma}

\newtheorem*{prop*}{Proposition}
\newtheorem{cor}[thm]{Corollary}

\theoremstyle{definition}

\newtheorem*{defn*}{Definition}
\newtheorem*{claim}{Claim}

\newtheorem*{ack}{Acknowledgement}
\theoremstyle{remark}
\newtheorem*{rem*}{Remark}
\def\A{\mathcal{A}}
\def\B{\mathcal{B}}
\def\D{\mathbf{\Delta}}
\def\S{\mathbf{\Sigma}}
\def\P{\mathbf{\Pi}}
\def\F{\mathcal{F}}
\def\I{\mathcal{I}}
\def\conc{\hat{\hbox to .4em{}}}
\def\l{\textnormal{length}}
\def\diam{\operatorname{diam}}
\def\osc{\operatorname{osc}}

\renewcommand{\labelenumi}{(\theenumi)}

\newtoks\by 
\newtoks\paper 
\newtoks\book 
\newtoks\jour 
\newtoks\yr 
\newtoks\no 
\newtoks\pages
\newtoks\vol 
\newtoks\publ 
\newtoks\publinfo
\newtoks\eds 
\newtoks\proc 
\newtoks\rada 
\newtoks\mathrev 
\newtoks\web
\def\ota{{\hbox{???}}} 
\def\cLear{\no=\ota\vol=\ota\publinfo=\ota\eds=\ota\proc=\ota\rada=\ota\web=\ota\by=\ota\paper=\ota\book=\ota\jour=\ota\yr=\ota 
\pages=\ota\vol=\ota\publ=\ota} 

\def\endpaper{\the\by, `\the\paper', 
\textit{\the\jour} \textbf{\the\vol} (\the\yr), no. {\the\no}, \the\pages.\cLear} 

\def\endpaperNN{\the\by, {\it \the\paper}, 
\the\jour{} {\bf \the\vol} (\the\yr),  \the\pages.\cLear}

\def\endbook{\the\by, {\em\the\book}, \the\publ, \the\yr.\cLear} 

\begin{document}

\maketitle

\begin{abstract}
We adapt a construction taken from `L. Motto Ros and B. Semmes, {\em A new proof 
of a theorem of Jayne and Rogers}, Real Anal. Exchange {\bf 35(1)} (2009/2010), 
195--204' in order to correct a mistake contained in the first part of the same paper. 
As a byproduct of the new construction, the Jayne-Rogers theorem is extended to functions whose range is a regular topological space, and a theorem of Solecki which sharpens the Jayne-Rogers theorem for separable metric spaces is extended to the non-separable context.
\end{abstract}

Let us first recall some important definitions together with the Jayne-Rogers theorem 
considered in \cite{mot-sem}.

\begin{defn*}
Let $X$ be a metric space and let $Y$ be a regular topological space. A function $f \colon X\to Y$ is said to be 
\emph{$\S^0_2$-measurable}, if $f^{-1}(U)\in\S^0_2(X)$ for every open set 
$U\subset Y$, and it is said to be a \emph{$\D^0_2$-function} if $f^{-1}(U)\in\D^0_2(X)$ for every open set $U\subset Y$
(equivalently, $f^{-1}(S)\in\S^0_2(X)$ for every $S\in\S^0_2(Y)$, if $Y$ is a metric space\footnote{In \cite{jaynerogers}, the \emph{first level Borel functions} were defined in this equivalent form \emph{for $Y$ a metric space}. However, if $Y$ is not metrizable (as it is the case for some results contained in the present work) then open sets in~$Y$ need not be $\S^0_2$, and hence there might be first level Borel functions  $f \colon X \to Y$ which are not piecewise continuous (even with $X$ metric). This shows that the change of the definition to $\D^0_2$-functions adopted in this paper provides a more natural notion in the non-metrizable context.}). 

The function $f$ is said to be \emph{piecewise continuous} if $X$ can be covered by 
a sequence $X_1,X_2,\ldots$ of closed sets such that $f\restriction X_n$ is continuous 
for every $n\in\omega$.

A function $f \colon X\to Y$ is said to be of 
\emph{Baire class~$1$} if it is the pointwise 
limit of a sequence of continuous functions $f_n \colon X\to Y$.

The metric space $X$ is said to be an \emph{absolute Souslin-$\F$} set if $X$ is a Souslin-$\F$ set in the completion $\widetilde{X}$ of $X$ under its metric, i.e., if it belongs to $\A\P^0_1(\widetilde{X})$, where $\A$ is the usual Souslin operation (see \cite[Definition~25.4]{kechris}). 
\end{defn*}

Notice that, by \cite[Theorem~25.7]{kechris}, separable absolute Souslin-$\F$ sets are exactly the \emph{analytic} sets.

\begin{thm}[Jayne-Rogers for $Y$ metric]\label{theorJR}
If $X$ is an absolute Souslin-$\F$ set and $Y$ is an arbitrary regular topological space, then $ f\colon X \to Y$ is a $\D^0_2$-function if and only if it is piecewise continuous.
\end{thm}

According to the authors of~\cite{jaynerogers}, their original proof of this theorem, with $Y$ assumed to be a metric space, is long and quite complicated (even in the separable case), but a more elementary and very short proof was presented some years ago in~\cite{mot-sem}. The proof was divided in two parts: in the first one, the general case of an absolute Souslin-$\F$ set $X$ was reduced to the particular situation in which $X$ is a zero-dimensional complete metric space, while in the second one (which contains the main construction of the paper), the Jayne-Rogers theorem was proved under such extra assumption on~$X$. 
However, as pointed out by the first author of this paper (and, independently, by M.~Sabok), the reduction used in the first part is faulty. To explain this in more detail, let us continue with some more definitions. 

Let $X$ be an absolute Souslin-$\F$ set and $Y$ be an arbitrary regular topological space. 
Given a function $f \colon X \to Y$, let us denote by $\I_f$ the $\sigma$-ideal of all subsets $A \subset X$ for which there is a set $S\in\S^0_2(X)$ such that $A \subset S$ and $f \restriction S$ is piecewise continuous. Note that this definition of $\I_f$ is equivalent to that given in~\cite{mot-sem} (up to the restriction to $X$ from its completion).

It is claimed on~\cite[p.~197]{mot-sem} that if $A \subset X$ is such that $A \notin \I_f$ then $f\restriction A$ is not piecewise continuous (where $f$ is implicitly assumed to be $\S^0_2$-measurable and $Y$ is a metric space). However, this is not true even if we further assume that $X$ be a Polish space and $A$ be a $\P^0_2(X)$ set, as the following counterexample shows. Consider $X = \R$, $A=\R \setminus \Q$, let $\langle q_n \mid n \in \N \rangle$ be any enumeration without repetitions of the rational numbers, and let $f \colon \R \to \R$ be the function
\[
f(x)=
\begin{cases}
\frac{1}{n}& \text{if } x=q_n \text{ for some }n \in \N,\\
0 & \text{if } x \in A.
\end{cases}
\]
Then $f$ is of Baire class $1$ (hence also $\S^0_2$-measurable), and by the Baire category theorem $A \notin \I_f$ since $f$ is discontinuous on every nonempty open set: but since $f$ is constant on~$A$, $f \restriction A$ is obviously piecewise continuous.

The mentioned reduction in the first part of~\cite{mot-sem} explicitly uses the claim above, so by the  previous counterexample  we must conclude that this part of the proof is wrong (even though the second part still works, so that the proof of the Jayne-Rogers theorem contained in~\cite{mot-sem} is still valid when $X$ is assumed to be a zero-dimensional complete metric space). 

A similar reduction had already been used at the beginning of the proof of Solecki's~\cite[Theorem~3.1]{sol}, a further sharpening of the Jayne-Rogers theorem in the context of \emph{separable} metric spaces. To precisely state this last result we must introduce one more notion: for $f \colon X_0 \to Y_0$ and $g \colon X_1 \to Y_1$, we say that \emph{$f$ is contained in $g$} ($f \sqsubseteq g$ in symbols) if and only if there are embeddings (i.e., open continuous injections)  $\phi \colon X_0 \to X_1$ and $\psi \colon f(X_0) \to Y_1$ such that 
\[ \psi \circ f = g \circ \phi. \] 

\begin{thm}[Solecki]\label{theorSolecki}
Let $X$ be an analytic  set, $Y$ be a separable metric space, and $f \colon X \to Y$ be a $\S^0_2$-measurable function. Then either $f$ is piecewise continuous or else one of $L, L_1$ is contained in $f$, where $L$ and $L_1$ are the two Lebesgue's functions defined on \cite[p.~522]{sol}.
\end{thm}

To prove this result, Solecki first observed on \cite[p.~530]{sol} that one can reduce the case of an arbitrary analytic set $X$ to the simpler case when $X$ is Polish, i.e., that if $f $ is not piecewise continuous then there is a $\P^0_2(\widetilde{X})$-set $G\subset X$ such that $f \restriction G$ is not piecewise continuous as well. However, in the proof presented in \cite{sol} such a set $G$ is obtained by implicitly using again the claim mentioned above, and by this reason it may fail to have the desired property (as the counterexample above shows, even in the separable case we have that $G \notin \I_f$ does not imply that $f \restriction G$ is not piecewise continuous). 

After obtaining the results contained in this note, in a recent conversation with Solecki we learned that the imprecision in \cite[Proof of Theorem~3.1]{sol} had already been pointed out to him by Sabok in 2010, and that he then straightaway found a very simple procedure which allows to manipulate the $\P^0_2(\widetilde{X})$-set $G \notin \I_f$ obtained in his original proof in order to get a new $\P^0_2(\widetilde{X})$-set $G^\prime \supset G$ for which $f \restriction G^\prime$ is not piecewise continuous, thus completely fixing the mentioned problem in \cite[Proof of Theorem~3.1]{sol} (for the sake of completeness, this argument is included in the remark at the end of this paper with prof. Solecki's kind permission). As later observed by Sabok, such an argument can be easily extended to the non-separable setting: however, such a generalization still does not fully fix the reduction used in~\cite[Proof of Theorem~1.1]{mot-sem} because the resulting set $G^\prime$ need not be zero-dimensional (so that we still do not have that $f \restriction G^\prime$ is of Baire class $1$, as required in \cite[Theorem~2.1]{mot-sem}).

The main goal of this paper is to completely fix the proof of Theorem~\ref{theorJR} contained in \cite{mot-sem}. However, the problematic reductions encountered in the proofs of both Theorem~\ref{theorJR} and Theorem~\ref{theorSolecki} also suggest that it would be interesting to abstractly know which kind of ``reductions'' can be used in those theorems, i.e., to ask under which extra conditions on $X$ and $Y$, a proof of one of them for such particular case automatically implies the general result.

Therefore, in what follows we will first fully reprove the Jayne-Rogers theorem~\ref{theorJR} by generalizing the construction presented in \cite[Proof of Theorem~2.1]{mot-sem} from the case of a zero-dimensional complete metric space $X$ to the case of an arbitrary absolute Souslin-$\F$ set. One benefit of this approach is that it allows us to weaken the assumption on~$Y$ from a metric space to a regular topological space. The other benefit is that, as a consequence of such generalization, we will get that not only reductions of the kind mentioned in this introduction are actually valid, but that in fact an even much stronger result of this type holds (see Theorem~\ref{theorreduction}): this provides an answer to the abstract question addressed at the end of the previous paragraph. Interestingly enough, combining this last result with the proof of \cite[Theorem~3.1]{sol} for Polish spaces, we will get as a corollary that Solecki's theorem~\ref{theorSolecki} holds even in the broader context of arbitrary absolute Souslin-$\F$ sets $X$ and arbitrary metric spaces $Y$ (see Theorem~\ref{theorSoleckigen}).

\medskip

We denote by $\omega$ the set of all non-negative integers. The set of all binary sequences of finite length is denoted by $2^{<\omega}$, while $2^\omega$ is the Cantor space. Similarly, $\omega^{<\omega}$ denotes the set of all sequences of non-negative integers of finite length and $\omega^\omega$ is the Baire space. From now on, unless said otherwise, we assume $X$ to be a metric space and $Y$ to be a regular topological space. For $A\subset X$, we denote the closure of $A$ in $X$ by $\overline{A}$ and the closure of~$A$ in~$\widetilde{X}$ by $\overline{A}^\ast$. For any undefined notation we refer to~\cite{kechris}.

Given $A,B\subset Y$, we say that $A$ and $B$ are \emph{strongly disjoint} if $\overline{A}\cap\overline{B}=\emptyset$. Let $f \colon X\to Y$ be a function. We put $A^f=f^{-1}(Y\setminus\overline{A})$. Note that one has $(A\cup B)^f=A^f\cap B^f$ and also note that if $A,B$ are strongly disjoint and $A^f,B^f\in\I_f$, then $X\in\I_f$, since $\{A^f,B^f\}$ is a covering of~$X$.

Let $x\in X$, $X^\prime\subset X$ and $A\subset Y$. We say that the pair $(x,X^\prime)$ is \emph{$f$-irreducible outside~$A$} if for every open neighborhood $V\subset X$ of~$x$ we have $A^f\cap X^\prime\cap V\not\in\I_f$. Otherwise we say that $(x,X^\prime)$ is \emph{$f$-reducible outside~$A$}. Notice that if $(x,X^\prime)$ is $f$-irreducible outside~$A$ then $x\in\overline{A^f\cap X^\prime}$, as otherwise $A^f\cap X^\prime\cap V=\emptyset\in\I_f$ for some open neighborhood $V$ of~$x$. Also notice that if $(x,X^\prime)$ is $f$-irreducible outside~$A$ then $(x,X^\prime\cap W)$ is $f$-irreducible outside~$A$ for any neighborhood $W$ of~$x$. Recall that a family $\B$ of subsets of~$X$ is said to be \emph{discrete} if $X$ can be covered by open sets each having a nonvoid intersection with at most one member of $\B$.

\begin{lem}
\label{lem1}
Suppose $f \colon X\to Y$ is a $\S^0_2$-measurable function, $X^\prime$ is a subset of~$X$ and $A\subset Y$ is an open set such that $X^\prime\subset A^f$. Then the following assertions are equivalent:
{
\renewcommand{\labelenumi}{\textnormal{(\roman{enumi})}}
\begin{enumerate}
\item
$X^\prime\not\in\I_f$;
\item
there is an $x\in\overline{X^\prime}$ and an open set $U\subset Y$ strongly disjoint from~$A$ such that $f(x)\in U$ and $(x,X^\prime)$ is $f$-irreducible outside~$U$.
\end{enumerate}
}
\end{lem}

\begin{proof}
(ii) $\Rightarrow$ (i): If $(x,X^\prime)$ is $f$-irreducible outside~$U$, then $U^f\cap X^\prime\cap X\not\in\I_f$ and therefore also $X^\prime\not\in\I_f$.

(i) $\Rightarrow$ (ii): Assume toward a contradiction that (ii) does not hold, i.e., for every $x\in\overline{X^\prime}$ and every open set $U\subset Y$ strongly disjoint from~$A$ such that $f(x)\in U$ we have that $(x,X^\prime)$ is $f$-reducible outside~$U$. Then there is some open neighborhood $V\subset X$ of~$x$ such that $U^f\cap X^\prime\cap V\in\I_f$. Let $\B$ be a $\sigma$-discrete base for the topology of~$X$, i.e., $\B=\bigcup_{n\in\omega}\B_n$ with each $\B_n$ discrete (see, e.g., \cite[Theorem~4.4.3]{eng}). Let $Q_n$ be the union of the elements $B$ of $\B_n$ such that $B\cap X^\prime\in\I_f$. Since, by~\cite[Lemma~2.2]{mot-sem}, $\I_f$ is closed under discrete unions, $Q_n\cap X^\prime\in\I_f$ for every $n\in\omega$. Finally, put $Q=\bigcup_n Q_n$ and notice that $Q\cap X^\prime\in\I_f$ and that $Q$ is open and contains as a subset each open set~$W$ for which $W\cap X^\prime\in\I_f$.

We claim that $f\restriction(\overline{X^\prime}\setminus Q)\cap A^f$ is continuous. Suppose otherwise, so that there is an $x\in(\overline{X^\prime}\setminus Q)\cap A^f$ and an open set $U\subset Y$ such that $f(x)\in U$ and there is no open neighborhood $V$ of~$x$ such that $f(V\cap(\overline{X^\prime}\setminus Q)\cap A^f)\subset U$. Using regularity of~$Y$, there is an open set $U^\prime\subset Y$ strongly disjoint from $A$ such that $f(x)\in U^\prime$ and $\overline{U^\prime}\subset U$. Let $V\subset X$ be an open neighborhood of~$x$ given by the failure of (ii) on the inputs $x$ and~$U^\prime$, so that $(U^\prime)^f\cap X^\prime\cap V\in\I_f$. By our hypothesis there is $x^\prime\in V\cap(\overline{X^\prime}\setminus Q)\cap A^f$ such that $f(x^\prime)\not\in \overline{U^\prime}$, and by using regularity of~$Y$ again, we can find an open neighborhood $U^{\prime\prime}\subset Y$ of~$f(x^\prime)$ strongly disjoint from $U^\prime$ and~$A$. Let now $V^\prime\subset X$ be an open set given by the failure of (ii) on inputs $x^\prime$ and~$U^{\prime\prime}$, so that $(U^{\prime\prime})^f\cap X^\prime\cap V^\prime\in\I_f$. By the strong disjointness of $U^\prime$ and $U^{\prime\prime}$, $\{(U^\prime)^f,(U^{\prime\prime})^f\}$ is a covering of the set $X^\prime\cap V\cap V^\prime$. Therefore $X^\prime\cap V\cap V^\prime\in\I_f$ and we must have $V\cap V^\prime\subset Q$. But this implies $x^\prime\in Q$, a contradiction!

Thus $f\restriction(\overline{X^\prime}\setminus Q)\cap A^f$ is continuous. Since $(\overline{X^\prime}\setminus Q)\cap A^f$ is a $\S^0_2(X)$ set containing~$X^\prime\setminus Q$, we get that $X^\prime\setminus Q\in\I_f$, and so $X^\prime\in\I_f$.
\end{proof}

\begin{lem}
\label{lem2}
Let $f \colon X\to Y$ be any function, $x\in X$, $X^\prime\subset X$, $A\subset Y$ and let $U_0,\ldots,U_n$ be a sequence of pairwise strongly disjoint open subsets of~$Y$. If $(x,X^\prime)$ is $f$-irreducible outside~$A$, then there is at most one $i\in\{0,\ldots,n\}$, such that $(x,X^\prime)$ is $f$-reducible outside~$A\cup U_i$.
\end{lem}
\begin{proof}
Assume that there are two indices $i,j\in\{0,\ldots,n\},\,i\neq j,$ such that $(x,X^\prime)$ is $f$-reducible outside both $A\cup U_i$ and $A\cup U_j$. Then there are open neighborhoods $V_i$ and $V_j$ of~$x$ such that $(A\cup U_i)^f\cap X^\prime\cap V_i\in\I_f$ and $(A\cup U_j)^f\cap X^\prime\cap V_j\in\I_f$. Since $U_i$ and $U_j$ are strongly disjoint, this implies that $A^f\cap X^\prime\cap V_i\cap V_j\in\I_f$, and thus $V_i\cap V_j$ contradicts the fact that $(x,X^\prime)$ is $f$-irreducible outside~$A$.
\end{proof}

\begin{lem}
\label{lem-main}
Let $X$ be an absolute Souslin-$\F$ set and $Y$ an arbitrary regular topological space. Let $f\colon X \to Y$ be a $\S^0_2$-measurable function which is not piecewise continuous. Then there is an open set $\widehat{U}\subset Y$ and a continuous reduction $g \colon 2^\omega\to X$ from $S=\{z\in 2^\omega \mid \exists i\forall j\ge i(z(j)=0)\}$ to $f^{-1}(\widehat{U})$.

Further, it is possible to construct this reduction so that it is in fact an embedding, $f \circ g$ is continuous on every $z \notin S$, and that for all $\varepsilon>0$, $\osc(f\circ g,z)< \varepsilon$ for all but finitely many $z \in S$.
\end{lem}

\begin{proof}
Since $X$ is an absolute Souslin-$\F$ set, we can write 
\[X = \bigcup_{\nu\in\omega^\omega}\bigcap_{n\in\omega}F_{\nu|n},\]
where $F_{\nu|n}$ is closed in the completion~$\widetilde{X}$ of~$X$ and $F_{\nu|n+1}\subset F_{\nu|n}$ for every $\nu\in\omega^\omega$ and $n\in\omega$. For every $\mu\in\omega^{<\omega}$ we denote
\[X_\mu = \bigcup_{\substack{\nu\in\omega^\omega\\\mu\subset\nu}}\bigcap_{n\in\omega}F_{\nu|n}.\]
Clearly, $X_\mu\subset F_\mu$ and $X_\mu=\bigcup_{n\in\omega}X_{\mu\conc n}$ for every $\mu\in\omega^{<\omega}$. The construction will be carried out by induction with respect to the order $\preceq$ on $2^{<\omega}$ defined by
\[s\preceq t \ \Longleftrightarrow\ \l(s)<\l(t) \vee (\l(s)=\l(t) \wedge s\le_{\textnormal{lex}}t),\]
where $\le_{\textnormal{lex}}$ is the usual lexicographical order on~$2^{\l(s)}$. We write $s\prec t$ if $s\preceq t$ and $s\neq t$.

We will construct a sequence $\langle V_s \mid s \in 2^{<\omega}\rangle$ of subsets of~$X$, a sequence $\langle x_s \mid s \in 2^{<\omega}\rangle$ of points of~$X$, a sequence $\langle U_s \mid s\in 2^{<\omega}\rangle$ of subsets of~$Y$, and a mapping $h \colon 2^{<\omega}\to \omega^{<\omega}$ such that for every $s\in 2^{<\omega}$:

\begin{enumerate}
\item
\label{c1}
if $t\subset s$ then $V_s\subset V_t$,
\item
\label{c2}
$\diam V_s\le 2^{-\l(s)+1}$,
\item
\label{c3}
$x_s\in\overline{V_s}$,
\item
\label{c4}
$f(x_s)\in U_s$,
\item
\label{c5}
$U_s$ is open,
\item
\label{c6}
if $s=t\conc 0$ then $x_s=x_t$,
\item
\label{c7}
if the last digit of $s$ is $1$ then $f(\overline{V_s})\cap\overline{\bigcup_{u\prec s}U_u}=\emptyset$,
\item
\label{c8}
$(x_t,V_t)$ is $f$-irreducible outside $A$ for every $t\preceq s$, where $A=\bigcup_{u\preceq s}U_u$,
\item
\label{c9}
if $t\subsetneq s$ then $h(t)\subset h(s)$ and if, moreover, the last digit of $s$ is $1$ then $h(t)\subsetneq h(s)$,
\item
\label{c10}
$V_s\subset X_{h(s)}$,
\item
\label{c11}
if the last digit of $s$ is $1$ then $\diam(f(\overline{V_s}))\le 2^{-\l(s)}$,
\item
\label{c12}
the family $\{V_t \mid t\in 2^n\}$ is pairwise strongly disjoint for every $n\in\omega$.
\end{enumerate}

At the first stage, let $x$ and $U$ be given as in Lemma~\ref{lem1} applied to $X^\prime=X$ and $A=\emptyset$. Then put $x_\emptyset=x$, $U_\emptyset=U$, $h(\emptyset)=\emptyset$ and let $V_\emptyset=B(x_\emptyset,1)$ be an open ball in~$X$ with the centre~$x_\emptyset$ and radius~$1$.

Suppose we have defined $V_t,x_t,U_t$ and $h(t)$ for every $t\prec s\conc 0$. Then putting $x_{s\conc 0}=x_s$, $U_{s\conc 0}=U_s$, $h(s\conc 0)=h(s)$, and $V_{s\conc 0}=V_s\cap B(x_s,2^{-\l(s)-1})$, it is easy to verify that conditions (\ref{c1})--(\ref{c11}) are satisfied.

Now suppose we have defined $V_t,x_t,U_t$ and $h(t)$ for every $t\prec s\conc 1$. Let $A=\bigcup_{t\prec s\conc 1}U_t$ and $O=Y\setminus\overline{A}$. By the inductive hypothesis, condition (\ref{c8}) applied with $t$ and $s$ replaced by $s$ and $s \conc 0$, respectively, says that $(x_s,V_s)$ is $f$-irreducible outside~$A$, so that in particular $f^{-1}(O)\cap V_s = A^f \cap V_s \notin \I_f$. Let\footnote{Note that if $Y$ is separable then this part of the argument can be simplified by taking a proper countable open decomposition of the set~$O$. The preimages form a countable family of $\S^0_2(X)$-sets and so the existence of the required set $C$ is immediate.} \label{inductionstep} $\{O^n_\alpha \mid n\in\omega,\,\alpha\in\Lambda_n\}$ be a $\sigma$-discrete open decomposition of~$O$ with $\diam(O^n_\alpha)\le 2^{-\l(s)-1}$ for each $n$ and $\alpha$ (see, e.g., \cite[Theorem~4.4.3]{eng}). By \cite[Theorem~3]{jaynerogers}, the inverse image under~$f$ of this family of sets is a discretely $\sigma$-decomposable family in~$X$, and so by a rather obvious modification of the proof of \cite[Theorem~4]{jaynerogers}, there are some $n$ and $\alpha$ such that $f^{-1}(O^n_\alpha)\cap V_s\notin\I_f$. Since $f$ is $\S^0_2$-measurable, we have $f^{-1}(O^n_\alpha)\in\S^0_2(X)$, and thus there exists a set $C\subset f^{-1}(O^n_\alpha)$ closed in $X$ such that $C\cap V_s\not\in\I_f$ and $\diam(f(C))\le 2^{-\l(s)-1}$. Also, since $X_{h(s)}=\bigcup_{n\in\omega}X_{h(s)\conc n}$ and (\ref{c10}) holds, there is some $n\in\omega$ such that $X^\prime=C\cap V_s\cap X_{h(s)\conc n}\not\in\I_f$. Put $h(s\conc1)=h(s)\conc n$.

\begin{claim}
There are $x_{s\conc1}\in\overline{X^\prime}$ and $U_{s\conc1}\subset Y$ such that $f(x_{s\conc1})\in U_{s\conc1}$, $U_{s\conc1}$ is open and strongly disjoint from $A$, $(x_t,V_t)$ is $f$-irreducible outside $A\cup U_{s\conc1}$ for every $t\prec s\conc1$ and $(x_{s\conc1},X^\prime)$ is $f$-irreducible outside $A\cup U_{s\conc1}$.
\end{claim}

\begin{proof}[Proof of the Claim.]
Let $k=|\{t\in 2^{<\omega} \mid t\prec s\conc1\}|$. Using Lemma~\ref{lem1}, for $j=0,\ldots,k$ recursively construct $x_j$ and $U_j$ so that $f(x_j)\in U_j$, $U_j$ is strongly disjoint from $A\cup U_{<j}$ (where $U_{<j}=\emptyset$ if $j=0$ and $U_{<j}=\bigcup_{i<j}U_i$ otherwise), and $(x_j,X^\prime\cap(U_{<j})^f)$ is $f$-irreducible outside $U_j$. Now notice that there must be some $\bar{\jmath}\in\{0,\ldots,k\}$, such that the claim is satisfied with $x_{s\conc1}=x_{\bar{\jmath}}$ and $U_{s\conc1}=U_{\bar{\jmath}}$: if not, by the pigeonhole principle there should be $j\neq j^\prime\le k$ and $t\prec s\conc1$ such that $(x_t,V_t)$ is $f$-reducible outside both $A\cup U_j$ and $A\cup U_{j^\prime}$, contradicting Lemma~\ref{lem2}.
\end{proof}

Finally, let $V_{s \conc 1}=X^\prime\cap B(x_{s\conc1},2^{-\l(s)-1})$. It is easy to check that all the conditions (\ref{c1})--(\ref{c11}) are satisfied. 

To achieve condition (\ref{c12}), we redefine the sets $\{V_t \mid t\in 2^n\}$ at the end of the induction step for $s=(1,1,\ldots,1)$ with $\l(s)=n$: It is clear from conditions (\ref{c3}), (\ref{c4}) and (\ref{c7}) that if $t,u\in 2^n$ and $t\neq u$ then $x_t\neq x_u$. Thus there is a family of nonempty pairwise strongly disjoint open balls $\{B_t \mid t\in 2^n\}$ centered in the respective $x_t$'s. By replacing $V_t$ by $V_t\cap B_t$ for every $t\in 2^n$, we get condition (\ref{c12}) and none of the other conditions is violated. This completes the recursive definition of the sequences required.

Now put $\widehat{U}=\bigcup_{s\in 2^{<\omega}}U_s$ and let $g \colon 2^\omega\to\widetilde{X}$ be defined by $g(z)=\bigcap_{n\in\omega}\overline{V_{z|n}}^\ast$ for every $z\in 2^\omega$. It follows from the completeness of~$\widetilde{X}$ and conditions (\ref{c1})--(\ref{c3}) that $g$ is well-defined and continuous. It remains to show that $g$ is an $X$-valued reduction from $S$ to $f^{-1}(\widehat{U})$.

If $z\in S$ then, by (\ref{c6}), for some $\bar{n}\in\omega$ we have that $x_{z|n}=x_{z|\bar{n}}=\bar{x}\in X$ for every $n\ge\bar{n}$. Therefore, by (\ref{c3}), $g(z)=\bar{x}$ and, by (\ref{c4}), $f(g(z))=f(\bar{x})\in U_{z|\bar{n}}\subset\widehat{U}$.

Assume now $z\not\in S$. By condition (\ref{c10}), $g(z)\in\overline{V_{z|n}}^\ast\subset F_{h(z|n)}$ for every $n\in\omega$. Therefore, using (\ref{c9}), there is some $\nu\in\omega^\omega$ such that $g(z)\in\bigcap_{n\in\omega}\overline{V_{z|n}}^\ast\subset\bigcap_{n\in\omega} F_{h(z|n)}= \bigcap_{m\in\omega} F_{\nu|m}\subset X$. Finally, since we have $g(z)\in\bigcap_{n\in\omega}\overline{V_{z|n}}^\ast\cap X=\bigcap_{n\in\omega}\overline{V_{z|n}}$, condition (\ref{c7}) implies that $f(g(z))\not\in\widehat{U}$.

Finally, condition (\ref{c11}) gives that $f \circ g$ is continuous on every $z \notin S$, and that for all $\varepsilon>0$, $\osc(f\circ g,z)< \varepsilon$ for all but finitely many $z \in S$. By (\ref{c12}), $g$ is an embedding.
\end{proof}

\begin{proof}[Proof of Theorem \ref{theorJR}.]
One direction is trivial. For the other direction, assume that $f$ is $\S^0_2$-measurable but not piecewise continuous. By Lemma~\ref{lem-main}, there is an open set $\widehat{U}\subset Y$ and a continuous reduction $g \colon 2^\omega\to X$ from $S=\{z\in 2^\omega \mid \exists i\forall j\ge i(z(j)=0)\}$ to $f^{-1}(\widehat{U})$. Since $S$ is a $\S^0_2$-complete set, $f^{-1}(\widehat{U})$ is also a $\S^0_2$-complete set and therefore not in~$\D^0_2(X)$.
\end{proof}

\begin{rem*}
The simplified version of the proof of Lemma~\ref{lem-main} without conditions (\ref{c11}) and (\ref{c12}) suffices to prove Theorem~\ref{theorJR}: in fact, these conditions are only used to ensure the additional statement of the lemma, which will be used only later in the alternative proof of Theorem \ref{theorSoleckigen}. Notice also that for this simplified version one can avoid any use of the results from \cite{jaynerogers} when considering the induction step for $s \conc 1$ of the main construction (see page \pageref{inductionstep}), because in this case it is enough to deduce from $A^f \cap V_s \notin \I_f$ and $A^f \in \S^0_2(X)$ (which follows from the $\S^0_2$-measurability of $f$) that there is a set $C \subset A^f$ closed in $X$ and such that $C \cap V_s \notin \I_f$.
\end{rem*}

The following result is a corollary of Lemma~\ref{lem-main}.

\begin{cor}\label{cor}
Let $X$ be an absolute Souslin-$\F$ set and $Y$ a regular topological space. If $f \colon X\to Y$ is $\S^0_2$-measurable and not piecewise continuous, then there is a Cantor set (i.e., a homeomorphic copy of the Cantor space $2^\omega$) $K \subset X$ such that $f \restriction K$ has the same properties. 
\end{cor}

\begin{proof}
Let $g$ be as in Lemma~\ref{lem-main}. Then $g(2^\omega)$ is a nonempty compact metrizable (hence Polish) zero-dimensional space such that $f \restriction g(2^\omega)$ is $\S^0_2$-measurable but not piecewise continuous (as $f^{-1}(\widehat{U})$ is a proper $\S^0_2$-set in $g(2^\omega)$). This in particular implies that $g(2^\omega)$ is uncountable. Using the Cantor-Bendixson theorem (see \cite[Theorem~6.4]{kechris}), write $g(2^\omega)$ as $P \cup C$, where $P$ is a (nonempty) perfect subset of $g(2^\omega)$ and $C$ is countable: then $P$ is a Cantor set by Brouwer's theorem (see \cite[Theorem 7.4]{kechris}) and $f \restriction P$ is $\S^0_2$-measurable but not piecewise continuous.
\end{proof}

\begin{thm}\label{theorreduction}
Let $X$ be an absolute Souslin-$\F$ set, $Y$ be a metric space, and $f \colon X \to Y$ be a $\S^0_2$-measurable function. If $f$ is not piecewise continuous then there are a Polish space $Y^\prime$ and a Baire class $1$ function $g \colon 2^\omega \to Y^\prime$ such that $g$ is not piecewise continuous and $g \sqsubseteq f$.
\end{thm}

\begin{proof}
By Corollary \ref{cor}, there exists a Cantor set $K \subset X$ such that 
$f^\prime =  f \restriction K$ is not piecewise continuous. 
By~\cite[Theorem~1]{frolik} and the 
subsequent proposition, $f^\prime(K)$ is Lindel\"of, and therefore a separable 
space (since $Y$ is a 
metric space). Thus, there is a Polish $Y^\prime \subset \widetilde{Y}$ which 
contains $f^\prime(K)$, and $f^\prime \colon K \to Y^\prime$ is 
a Baire class $1$ function because, being $K$ and $Y^\prime$ separable metrizable with $K$ zero-dimensional, 
this notion coincides with $\S^0_2$-measurability by classical results (see e.g.~\cite[Theorem 24.10]{kechris})\footnote{Alternatively, instead of using the result from \cite{frolik} one 
could also first observe by using~\cite[Theorem 8]{hansell} that $f \colon K \to Y$ is a Baire class $1$ function, and then directly show by induction on \( \alpha < \omega_1 \) that for every 
Baire class $\alpha$  function $f$  from a compact topological space $K$ to a 
metrizable space $Y$ there is a separable and closed (hence Polish) subspace 
$Y^\prime$ of the completion~$\widetilde{Y}$ of~$Y$ such that 
$f(K) \subset Y^\prime$ and $f \colon K \to Y^\prime$ is still of Baire class 
$\alpha$ (in particular, $f(K)$ is a separable subspace of $Y$).}. Let $h \colon 2^\omega \to K$ be any homeomorphism, and put 
$g = f^\prime \circ h \colon 2^\omega \to Y^\prime$. Then $g$ cannot be 
piecewise continuous (otherwise $f^\prime$ would be piecewise continuous as 
well) and $g \sqsubseteq f$, as witnessed by $h$ and the identity function on 
$g(2^\omega) = f^\prime(K)$.
\end{proof}

In particular, we get the following strengthening of Solecki's Theorem \ref{theorSolecki}.

\begin{thm}\label{theorSoleckigen}
Let $X$ be an absolute Souslin-$\F$ set, $Y$ be a metrizable space, and $f$ be a $\S^0_2$-measurable function. Then exactly one of the following holds:
{
\renewcommand{\labelenumi}{\textnormal{(\roman{enumi})}}
\begin{enumerate}
\item
$f$ is piecewise continuous;
\item
$L \sqsubseteq f$ or $L_1 \sqsubseteq f$.
\end{enumerate}
}
\end{thm}

\begin{proof}
As a consequence of Theorem \ref{theorreduction}, we get that the proof can be reduced to the case of a Baire class $1$ function $g \colon 2^\omega \to Y^\prime$ with $Y^\prime$ a Polish space. To see this, observe that:
\begin{enumerate}
\item If $g \sqsubseteq f$ and $g$ is not a $\D^0_2$-function then $f$ is not a $\D^0_2$-function;
\item the relation $\sqsubseteq$ is transitive.
\end{enumerate}
Now, it is enough to apply the argument contained in~\cite[Proof of Theorem 3.1]{sol} (for the special case of a Polish space $X$) to the function $g$ given by Theorem \ref{theorreduction} --- this is possible as $g$ is defined on the Cantor space $2^\omega$ and $Y^\prime$ is metric separable.
\end{proof}

\begin{proof}[Alternative proof]
Using the reduction $g$ given by Lemma~\ref{lem-main}, we have that $f \circ g$ is continuous on every $z \notin S$, and that for all $\varepsilon>0$, $\osc(f\circ g,z)< \varepsilon$ for all but finitely many $z \in S$. So we can now use the last thirteen lines of Solecki's \cite[Proof of Theorem~3.1]{sol} (which use only elementary or well-known results) to show that one of $L,L_1$ embeds into $f \circ g$, which implies the same for $f$ because $g$ is an embedding (by Lemma \ref{lem-main} again).
\end{proof}

\begin{rem*}

Let $X$ be an analytic set, $Y$ be a separable metric space and $f \colon X \to Y$ be a $\S^0_2$-measurable function which is not piecewise continuous. By~\cite[Theorem~1]{sol2}, there is a $\P^0_2(\widetilde{X})$-set $G\subset X$ such that $G \notin \I_f$. As already discussed, it is essentially  claimed on~\cite[p.~530]{sol} that $f \restriction G$ is not piecewise continuous, but this does not directly follow from $G \notin \I_f$. The following simple argument, which fixes this problem, is due to Solecki.

Without loss of generality, we can assume that no non-empty open subset of $G$ is in $\I_f$. We represent $G$ as the intersection of a decreasing sequence $\langle U_n \mid n \in \N \rangle$ of open subsets of~$\widetilde{X}$. Let $\langle V_n \mid n \in \N \rangle$ be an enumeration of a topological basis of $G$ consisting of non-empty sets. From the assumption on~$G$, we have that $f$ is discontinuous on 
$\overline{V}_n\cap U_n$ for every $n \in \N$. Let  $K_n\subset \overline{V}_n\cap U_n$ be a compact set on which $f$ is discontinuous (i.e., $K_n$ is some convergent sequence together with its limit). Then 
$G^\prime  = G\cup \bigcup_n K_n$  
is still a $\P^0_2(\widetilde{X})$-set since 
$G\cup \bigcup_n K_n  = \bigcap_n (U_n\cup \bigcup_{m<n} K_m)$.
Discontinuity points of $f$ restricted to $G^\prime$ are dense in it, so $f \restriction G^\prime$ is not piecewise continuous by the Baire category theorem.

As already discussed in the introduction, despite the fact that it can be easily extended to the non-separable setting, the argument considered in this remark does not fully fix the reduction used in~\cite[Proof of Theorem 1.1]{mot-sem}: this is why we think that it was nonetheless worth carrying out the proof contained in this paper. Moreover, we underline that, in any case, the approach we took in the present note has some advantages with respect to the kind of modifications considered in this remark, namely: a) it gives an even simpler, more general and more elementary proof of Jayne-Rogers theorem \ref{theorJR}, as it avoids any reference to nontrivial results like~\cite[Theorem~1]{sol2} or its non-separable version~\cite[Theorem 1.3]{holzajzel}, and b) it allows to extend Solecki's theorem~\ref{theorSolecki} to the non-separable context.
\end{rem*}

\begin{ack}
The authors would like to thank the referee for helpful suggestions on the preliminary version of this paper.
\end{ack}

\end{document}